\documentclass[12pt,twoside]{amsart}
\usepackage{amssymb}

\usepackage[all]{xy}
\nonstopmode

\numberwithin{equation}{section} 

\font\tengothic=eufm10 scaled\magstep 1 \font\sevengothic=eufm7
scaled\magstep 1
\newfam\gothicfam
       \textfont\gothicfam=\tengothic
       \scriptfont\gothicfam=\sevengothic

\newtheorem{theorem}{Theorem}[section]
\newtheorem{lemma}[theorem]{Lemma}
\newtheorem{proposition}[theorem]{Proposition}

\newtheorem{conjecture}[theorem]{Conjecture}

\theoremstyle{definition}
\newtheorem{definition}[theorem]{Definition} 

\newtheorem{example}[theorem]{Example}

\newcommand{\cD}{{\mathcal D}}
\newcommand{\cA}{{\mathcal A}}
\newcommand{\cB}{{\mathcal B}}

\newcommand{\cU}{{\mathcal U}}

\newcommand{\cC}{{\mathcal C}}

\newcommand {\ZZ}{\mathbb{Z}}

\begin{document}
\title[]{A note on the multiplicity of determinantal ideals}

\author[ Rosa M.\ Mir\'o-Roig]{ Rosa M.\
Mir\'o-Roig$^{*}$}

\address{Facultat de Matem\`atiques,
Departament d'Algebra i Geometria, Gran Via de les Corts Catalanes
585, 08007 Barcelona, SPAIN } \email{miro@ub.edu}

\date{\today}
\thanks{$^*$ Partially supported by MTM2004-00666.}



\begin{abstract} Herzog, Huneke, and Srinivasan have conjectured that for
any homogeneous $k$-algebra, the multiplicity
 is bounded above by a  function of the maximal degrees of the syzygies and below
 by a  function of the minimal degrees of the syzygies.
The goal of this paper is to establish the multiplicity conjecture
of Herzog, Huneke, and Srinivasan about the multiplicity of graded
Cohen-Macaulay algebras over a field $k$ for $k$-algebras $k[x_1,
\cdots ,x_n]/I$ being $I$ a  determinantal ideal of arbitrary
codimension.

\end{abstract}


\maketitle

\tableofcontents


  \section{Introduction} \label{intro}
Let $R=k[x_1, \cdots ,x_n]$ be a polynomial ring in $n$ variables
over a filed $k$, let $deg(x_i)=1$ and let $I\subset R$ be a
graded ideal of arbitrary codimension. Consider the minimal graded
free $R$-resolution of $R/I$:

$$ 0 \longrightarrow  \oplus _{j\in \ZZ}R(-j)^{\beta_{p,j}(R/I)}\longrightarrow
\cdots  \longrightarrow   \oplus _{j\in
\ZZ}R(-j)^{\beta_{1,j}(R/I)} \longrightarrow R\longrightarrow R/I
\longrightarrow 0$$ where we denote
$\beta_{i,j}(R/I)=Tor_{i}^R(R/I,k)_{j}$ the graded Betti number of
$R/I$. Many important numerical invariants of $I$ and the
associated scheme can be read off from the minimal graded free
$R$-resolution of $R/I$. For instance,  the Hilbert polynomial,
and hence the multiplicity $e(R/I)$ of $I$, can be written down in
terms of the shifts $j$ such that $\beta_{i,j}(R/I)\ne 0$ for some
$i$, $1\le i \le p$.

Let $c$ denote the codimension of $R/I$. Then $c\le p$ and
equality holds if and only if $R/I$ is Cohen-Macaulay. We define
$m_i(I)=\text{min} \{ j\in \ZZ \mid \beta_{i,j}(R/I)\ne 0 \}$ the
minimum degree shift at the $i$-th  step and $M_i(I)=\text{max} \{
j\in \ZZ \mid \beta_{i,j}(R/I)\ne 0 \}$ the maximum degree shift
at the $i$-th  step. We will simply write $m_i$ and $M_i$ when
there is no confusion. If $R/I$ is Cohen-Macaulay and has a pure
resolution, i.e. $m_i=M_i$ for all $i$, $1\le i \le c$, then
Huneke and Miller showed in \cite{HM} that
$$ e(R/I)=\frac{\prod_{i=1}^cm_i}{c!}.$$

Generalizing their result Herzog, Huneke, and Srinivasan made the
following multiplicity conjecture:

\begin{conjecture}\label{conjetura} If $R/I$ is Cohen-Macaulay then
$$\frac{\prod_{i=1}^cm_i}{c!}\le e(R/I)\le
\frac{\prod_{i=1}^cM_i}{c!}.$$
\end{conjecture}

Conjecture \ref{conjetura} has been extensively studied, and
partial results have been obtained. It turns out to be true for
the following type of ideals:
\begin{itemize}
\item Complete intersections \cite{HS} \item Powers of complete
intersection ideals \cite{GV} \item Perfect ideals with a pure
resolution \cite{HM} \item Perfect ideals with a quasi-pure
resolution (i.e. $m_i\ge M_{i-1}$) \cite{HS} \item Perfect ideals
 of codimension 2 \cite{HS} \item Gorenstein ideals of codimension 3 \cite{MNR}
 \item Perfect stable monomial ideals
 \cite{HS} \item Perfect square free strongly stable monomial
 ideals \cite{HS}.
\end{itemize}

 \vskip 2mm The goal of this paper is to prove Conjecture \ref{conjetura}
 for determinantal ideals of arbitrary codimension $c$, i.e.
 ideals generated by the maximal minors of a $t\times (t+c-1)$ homogeneous
polynomial matrix.  Determinantal ideals have been a central topic
in both commutative algebra and algebraic geometry and, due to
their important role, their study has attracted many researchers
and has received considerable attention in the literature.  Some
of the most remarkable results about determinantal ideals are due
to J.A. Eagon and M. Hochster in \cite{e-h}, and to J.A. Eagon and
D.G. Northcott in \cite{e-n}. J.A. Eagon and M. Hochster proved
that generic determinantal ideals are perfect. J.A. Eagon and D.G.
Northcott constructed a finite graded free resolution for any
determinantal ideal and, as a corollary, they got that
determinantal ideals are perfect. Since then many authors have
made important contributions to the study of determinantal ideals
and the reader can look at \cite{b-v},  \cite{BH}, \cite{KM} and
\cite{eise} for background, history  and a list of important
papers.

\vskip 2mm In this short note we verify that determinantal ideals
$I$ satisfy Herzog-Huneke-Srinivasan Conjecture which relates the
multiplicity $e(R/I)$ to the minimal and maximal shifts in the
graded minimal $R$-resolution of $R/I$.

\vskip 2mm \vskip 2mm Next we outline the structure of the paper.
In section 2, we first recall the basic facts on determinantal
ideals $I$
 of codimension $c$ defined by the maximal
minors of a $t\times (t+c-1)$ homogeneous matrix $\cA$ and the
associated complexes needed later on. We determine the minimal and
maximal shifts in the graded minimal free $R$-resolution of $R/I$
in terms of the degree matrix $\cU$ of $\cA$ and we state some
technical lemmas used in the inductive process of the proof of our
main Theorem (cf Theorem \ref{maintheorem}).

\vskip 2mm Section 3 is completely devoted to proving Conjecture
\ref{conjetura} for determinantal ideals $I$ of arbitrary
codimension. To prove it we use induction on the codimension $c$
of $I$ and for any $c$ induction on the size $t$ of the
homogeneous $t\times (t+c-1)$ matrix whose maximal minors generate
$I$ by successively deleting columns and rows of the largest
possible degree when we prove the lower bound and columns and rows
of the smallest possible degree when we prove the upper bound. We
end the paper with an example which illustrates that the upper and
lower bounds for the multiplicity $e(R/I)$ of a determinantal
ideal $I$ given in Theorem \ref{maintheorem} are sharp.

 \vskip 4mm {\bf Acknowledgement} The author thanks Laura Costa
 for all her help.

 \vskip 4mm

\section{Determinantal ideals}

In the first part of this  section, we  provide the background and
basic results on determinantal ideals needed in the sequel, and we
refer to \cite{b-v} and \cite{eise} for more details.

\vskip 2mm

Let $\cA$ be a homogeneous matrix, i.e. a matrix representing a
degree 0 morphism $\phi :F \longrightarrow G$ of free graded
$R$-modules. In this case, we denote by $I(\cA)$ the ideal of $R$
generated by the maximal minors of $\cA$.

\begin{definition}
 An homogeneous ideal $I\subset R$ of codimension $c$ is
called a \emph{determinantal} ideal if $I=I(\cA)$ for some
$t\times (t+c-1)$ homogeneous matrix $\cA$.
 \end{definition}

\vskip 2mm
 Let $I\subset R$  be a determinantal ideal
 of codimension    $c$
generated by the  maximal minors of a $t\times (t+c-1)$ matrix
$\cA=(f_{ji})_{i=1,...,t}^{j=1,...,t+c-1}$ where $f_{ji}\in {
k}[x_{1},...,x_{n}]$ are homogeneous polynomials of degree
$a_j-b_{i}$. We assume without loss of generality that $\cA$ is
minimal; i.e., $f_{ji}=0$ for all $i,j$ with $b_{i}=a_{j}$. If we
let $u_{j,i}=a_j-b_i$ for all $j=1, \dots , t+c-1$ and $i=1, \dots
, t$, the matrix $\cU=(u_{j,i})_{i=1,...t}^{j=1,...,t+c-1}$ is
called the {\em degree matrix} associated to $I$. By  re-ordering
degrees, if necessary,   we may also assume that $b_1 \ge ... \ge
b_t$ and $a_0 \le a_1\le ... \le a_{t+c-2}$. In particular, we
have:

\begin{equation} \label{order} u_{j,i}\le u_{j+1,i} \quad \text{ and }
\quad  u_{j,i}\le u_{j,i+1} \text{ for all } i,j.
\end{equation}

 \vskip 2mm  Note that the
degree matrix $\cU$ is completely determined by  $u_{1,1}$,
$u_{2,1}$, ... , $u_{c,1}$,   $u_{2,2}$, $u_{3,2}$, ... ,
$u_{c+1,2}$, ..., $u_{t,t}$, $u_{t+1,t}$, ... , $u_{c+t-1,t}$.
Moreover, the graded Betti numbers in the minimal free
$R$-resolution of
 $R/I(\cA)$ depend only upon the integers
 $$\{ u_{j,i} \} _{1\le i \le t}^{i\le j \le c+i-1}\subset \{u_{j,i}\}_{i=1,...t}^{j=1,...,t+c-1}$$ as described
 below.

\vskip 4mm
\begin{proposition}\label{miMi}
 Let $I\subset R$ be a determinantal ideal of codimension $c$ with
 degree matrix $\cU=(u_{ji})_{i=1,...t}^{j=1,...,t+c-1}$ as above.
 Then we have:
 \vskip 2mm
 \begin{itemize}
 \item[(1)] $m_i=u_{1,1}+u_{2,1}+\cdots +u_{i,1}+u_{i+1,2}+u_{i+2,3}+\cdots +u_{t+i-1,t} $ for  $1\le i \le c$,
 \item[(2)] $M_i=  u_{c-i+1,1}+u_{c-i+2,2}+\cdots +u_{t+c-i,t}+u_{t+c-i+1,t}+u_{t+c-i+2,t}+\cdots +u_{t+c-1,t} $ for  $1\le i \le c$.
 \end{itemize}
 \end{proposition}

\begin{proof}  We denote by
$\varphi :F \longrightarrow G$ the morphism of free graded
$R$-modules of rank $t$ and $t+c-1$, defined by the homogeneous
matrix $\cA$ associated to $I$. The Eagon-Northcott complex
${\cD}_0(\varphi^*):$

$$0 \longrightarrow \wedge^{t+c-1}G^* \otimes S_{c-1}(F)\otimes \wedge^tF\longrightarrow
\wedge^{t+c-2} G ^*\otimes S _{c-2}(F)\otimes \wedge
^tF\longrightarrow \ldots \longrightarrow$$ $$\wedge^{t}G^*
\otimes S_{0}(F)\otimes \wedge^tF\longrightarrow R \longrightarrow
R \longrightarrow R/I \longrightarrow 0 $$

\vskip 2mm \noindent gives us a graded minimal free
$R$-resolution of $R/I$ (See, for instance \cite{b-v}; Theorem
2.20 and \cite{eise}; Corollary A2.12 and Corollary A2.13). Now
the result follows after an straightforward computation.
\end{proof}

\vskip 2mm We will now fix the notation and prove the technical
lemmas needed  in the induction process we will use in next
section for proving the multiplicity Conjecture for determinantal
ideals of arbitrary codimension.

\vskip 2mm Let $I\subset R$ be a homogeneous ideal of codimension
$c$. Assume that $I$ is determinantal and let $\cA$ (resp $\cU$)
be the $t\times (t+c-1)$ homogeneous matrix (resp. degree matrix)
associated to $I$. Let $\cA'$ (resp $\cU'$) be the $(t-1)\times
(t+c-2)$ homogeneous matrix (resp. degree matrix) obtained
deleting the last column and the last row of $\cA$ and denote by
$I'$ the codimension $c$ determinantal ideal generated by the
maximal minors of $\cA'$. Since the multiplicity of $R/I$ and
$R/I'$ are completely determined by the corresponding degree
matrices, it is enough to consider an example of ideal for any
degree matrix. So, from now on, we take

$$\cA:= \begin{pmatrix}
x_1 ^{u_{1,1}}& x_2 ^{u_{2,1}}& \cdots & x_{c-1}^{u_{c-1,1}} & x_{c} ^{u_{c,1}}& 0 & 0 & \cdots & 0 & 0 \\
0 & x_1 ^{u_{2,2}}& x_2 ^{u_{3,2}}& \cdots & x_{c-1}^{u_{c,2}} & x_{c} ^{u_{c+1,2}}& 0 & \cdots & 0 & 0\\
0 & 0 & x_1 ^{u_{3,3}}& x_2 ^{u_{4,3}}&  \cdots & x_{c-1}^{u_{c+1,3}} & x_{c} ^{u_{c+2,3}}& \cdots & 0 & 0\\
\vdots &  \vdots & \vdots & \vdots &  \vdots & \vdots &  \vdots & \vdots & \vdots & \vdots\\
 0 & 0&
\cdots & 0 & x_1 ^{u_{t-1,t-1}}& x_2^{u_{t,t-1}} & \cdots & x_{c-1} ^{u_{c+t-3,t-1}}& x_{c}^{u_{c+t-2,t-1}} & 0\\
0 &
 0 &
0 & \cdots & 0 & x_1^{u_{t,t}} & x_2 ^{u_{t+1,t}}& \cdots &
x_{c-1}^{u_{t+c-2,t}} & x_{c}^{u_{t+c-1,t}}
\end{pmatrix}$$
and $$\cA':=\begin{pmatrix}x_1 ^{u_{1,1}}& x_2 ^{u_{2,1}}& \cdots & x_{c-1}^{u_{c-1,1}} & x_{c} ^{u_{c,1}}& 0 & 0 & \cdots & 0 \\
0 & x_1 ^{u_{2,2}}& x_2 ^{u_{3,2}}& \cdots & x_{c-1}^{u_{c,2}} & x_{c} ^{u_{c+1,2}}& 0 & \cdots & 0 \\
0 & 0 & x_1 ^{u_{3,3}}& x_2 ^{u_{4,3}}&  \cdots & x_{c-1}^{u_{c+1,3}} & x_{c} ^{u_{c+2,3}}& \cdots & 0 \\
\vdots &  \vdots & \vdots & \vdots &  \vdots & \vdots &  \vdots & \vdots & \vdots \\
 0 & 0&
\cdots & 0 & x_1 ^{u_{t-1,t-1}}& x_2^{u_{t,t-1}} & \cdots & x_{c-1} ^{u_{c+t-3,t-1}}& x_{c}^{u_{c+t-2,t-1}} \\
\end{pmatrix}$$
 \vskip 2mm Let $J\subset R$ be the codimension $c-1$
determinantal ideal generated by the maximal minors of the
$t\times (t+c-2)$ homogeneous matrix

$$\cB:=\begin{pmatrix}
x_1 ^{u_{1,1}}& x_2 ^{u_{2,1}}& \cdots & x_{c-1}^{u_{c-1,1}} & x_{c} ^{u_{c,1}}& 0 & 0 & \cdots & 0  \\
0 & x_1 ^{u_{2,2}}& x_2 ^{u_{3,2}}& \cdots & x_{c-1}^{u_{c,2}} & x_{c} ^{u_{c+1,2}}& 0 & \cdots & 0 \\
0 & 0 & x_1 ^{u_{3,3}}& x_2 ^{u_{4,3}}&  \cdots & x_{c-1}^{u_{c+1,3}} & x_{c} ^{u_{c+2,3}}& \cdots & 0 \\
\vdots &  \vdots & \vdots & \vdots &  \vdots & \vdots &  \vdots & \vdots & \vdots \\
 0 & 0&
\cdots & 0 & x_1 ^{u_{t-1,t-1}}& x_2^{u_{t,t-1}} & \cdots & x_{c-1} ^{u_{c+t-3,t-1}}& x_{c}^{u_{c+t-2,t-1}} \\
0 &
 0 &
0 & \cdots & 0 & x_1^{u_{t,t}} & x_2 ^{u_{t+1,t}}& \cdots &
x_{c-1}^{u_{t+c-2,t}}
\end{pmatrix}$$

\noindent obtained deleting the last column of $\cA$. Analogously,
we consider  $\cA''$ (resp $\cU''$)  the $(t-1)\times (t+c-2)$
homogeneous matrix (resp. degree matrix) obtained deleting the
first column and the first row of $\cA$ and we denote by $I''$ the
codimension $c$ determinantal ideal generated by the maximal
minors of $\cA''$. Let $\cC$ be the $t\times (t+c-2)$ homogeneous
matrix obtained deleting the first column  of $\cA$ and let
$K\subset R$ be the codimension $c-1$ determinantal ideal
generated by the maximal minors of $\cC$.

The ideal $I$ is obtained from $I'$ by a basic double G-link as
well as from $K$ by a basic double G-link. Indeed, we have

\vskip 2mm
\begin{lemma}\label{key} With the above notation, it holds
\begin{itemize}
\item[(1)] $I=J+x_c^{u_{t+c-1,t}}I'$ and $I=K+x_1^{u_{1,1}}I''$.
\item[(2)] The sequences
$$0\longrightarrow J(-u_{t+c-1,t})\longrightarrow
I'(-u_{t+c-1,t})\oplus J \longrightarrow J+
x_c^{u_{t+c-1,t}}I'=I\longrightarrow 0
$$ and $$0\longrightarrow K(-u_{1,1})\longrightarrow
I''(-u_{1,1})\oplus K \longrightarrow K+
x_1^{u_{1,1}}I''=I\longrightarrow 0
$$ are exact.\item[(3)] $e(R/I)=e(R/I')+u_{t+c-1,t}\cdot e(R/J)$ and
$e(R/I)=e(R/I'')+u_{1,1}\cdot e(R/K)$.
\end{itemize}
\end{lemma}
\begin{proof} (1) The equalities of ideals are immediate.

(2) and (3) follow from \cite{KMMNP}; Lemma 4.8.
\end{proof}

\begin{lemma}\label{keyinduction} With the above notation, we
have
\begin{itemize}
\item[(1)] $m_i=m_i(I)=m_i(I')+u_{t+i-1,t}=m_i'+u_{t+i-1,t}$ for
all $1\le i \le c$, \item[(2)]
$M_i=M_i(I)=M_i(I'')+u_{c-i+1,1}=M_i''+u_{c-i+1,1}$ for all $1\le
i \le c$,\item[(3)] $m_i(J)=m_i(I)=m_i$ for all $1\le i \le c-1$,
and \item[(4)] $M_i(K)=M_{i}(I)=M_i$ for all $1\le i \le c-1$.
\end{itemize}
\end{lemma}
\begin{proof} It follows from Proposition \ref{miMi}.
\end{proof}

 \vskip 2mm


\section{The multiplicity Conjecture}

Using the fact that the ideal $I$ is obtained from the ideal $I'$
(resp. $I''$) by a basic double G-link, we can now show that
Conjecture 1.1 is true for determinantal ideals of arbitrary
codimension.

\vskip 2mm
\begin{theorem}\label{maintheorem} Let $I\subset R$ be a
determinantal ideal of codimension $c$. Then the following lower
and upper bounds hold:
\begin{itemize}
\item[(1)] $e(R/I)\ge \frac{\prod_{i=1}^cm_i}{c!}$, and \item[(2)]
$e(R/I)\le \frac{\prod_{i=1}^cM_i}{c!}$.
\end{itemize}
\end{theorem}
\begin{proof} As we explained in section 2, it is enough to
 prove the result for the ideal $I$ generated by the maximal minors of the $t\times
 (t+c-1)$ matrix

 $$\cA:= \begin{pmatrix}
x_1 ^{u_{1,1}}& x_2 ^{u_{2,1}}& \cdots & x_{c-1}^{u_{c-1,1}} & x_{c} ^{u_{c,1}}& 0 & 0 & \cdots & 0 & 0 \\
0 & x_1 ^{u_{2,2}}& x_2 ^{u_{3,2}}& \cdots & x_{c-1}^{u_{c,2}} & x_{c} ^{u_{c+1,2}}& 0 & \cdots & 0 & 0\\
0 & 0 & x_1 ^{u_{3,3}}& x_2 ^{u_{4,3}}&  \cdots & x_{c-1}^{u_{c+1,3}} & x_{c} ^{u_{c+2,3}}& \cdots & 0 & 0\\
\vdots &  \vdots & \vdots & \vdots &  \vdots & \vdots &  \vdots & \vdots & \vdots & \vdots\\
 0 & 0&
\cdots & 0 & x_1 ^{u_{t-1,t-1}}& x_2^{u_{t,t-1}} & \cdots & x_{c-1} ^{u_{c+t-3,t-1}}& x_{c}^{u_{c+t-2,t-1}} & 0\\
0 &
 0 &
0 & \cdots & 0 & x_1^{u_{t,t}} & x_2 ^{u_{t+1,t}}& \cdots &
x_{c-1}^{u_{t+c-2,t}} & x_{c}^{u_{t+c-1,t}}
\end{pmatrix}$$

 \vskip 2mm  (1) We proceed by induction on the codimension $c$ of
$I$. If $c=1$ then $I$ is a principal ideal  and the result is
trivial. For $c=2$ the result was proved by  Herzog and Srinivasan
in \cite{HS}. Assume $c\ge 3$. We will now induct on $t$. If $t=1$
then $I$ is a complete intersection ideal and hence the result is
well known. Assume $t>1$. Let $\cA'$ (resp. $\cB$) be the matrix
obtained deleting the last column and the last row (resp. the last
column) of the matrix $\cA$ and let $I'$ (resp. $J$) be the ideal
generated by the maximal minors of $\cA'$ (resp. $\cB$). Let
$m_i$, $m_i'$  and $m_i(J)$ be the minimal shifts in the graded
minimal free $R$-resolution of $R/I$, $R/I'$ and $R/J$,
respectively (see Proposition \ref{miMi} and Lemma
\ref{keyinduction}).

\vskip 2mm By Lemma \ref{key} (3),
$$e(R/I)=e(R/I')+u_{t+c-1,t}\cdot e(R/J),$$ by hypothesis of
induction on $c$ and Lemma \ref{keyinduction} (3), we have
$$e(R/J)\ge
\frac{\prod_{i=1}^{c-1}m_i(J)}{(c-1)!}=\frac{\prod_{i=1}^{c-1}m_i}{(c-1)!},$$
and by hypothesis of induction on $t$ we have $$e(R/I')\ge
\frac{\prod_{i=1}^cm_i'}{(c)!}.$$

Therefore, since $m_i=m_i'+u_{t+i-1,t}$ (Lemma \ref{keyinduction}
(1)), we have $$c!e(R/I)\ge \prod_{i=1}^cm_i$$ if and only if
$$\prod_{i=1}^cm_i'+cu_{t+c-1,t}\prod_{i=1}^{c-1}m_i\ge
\prod_{i=1}^cm_i=$$
$$\prod_{i=1}^c(m_i'+u_{t+i-1,t})=$$ $$u_{t+c-1,t}\prod_{i=1}^{c-1}m_i+
\prod_{i=1}^cm_i'+m_c'\sum_{r=0}^{c-2}(u_{t+r,t}m_1\cdots
m_rm_{r+2}'\cdots m_{c-1}')=$$ $$u_{t+c-1,t}\prod_{i=1}^{c-1}m_i+
\prod_{i=1}^cm_i'+\sum_{r=0}^{c-2}(u_{t+r,t}m_1\cdots
m_rm_{r+2}'\cdots m_{c-1}'m_c')
 $$ if and only if
$$(c-1)u_{t+c-1,t}\prod_{i=1}^{c-1}m_i\ge \sum_{r=0}^{c-2}(u_{t+r,t}m_1\cdots
m_rm_{r+2}'\cdots m_{c-1}'m_c').
$$

Since,  for all integer $ i$, $ 1\le i \le c-1$, and  for all
integer $ r$, $ 0\le r \le c-2$, we have the inequalities

$$m_{i}-m_{i+1}'=
(u_{1,1}+u_{2,1}+\cdots +u_{i,1}+u_{i+1,2}+u_{i+2,3}+\cdots
+u_{t+i-1,t})-$$
$$(u_{1,1}+u_{2,1}+\cdots +u_{i,1}+u_{i+1,1}+u_{i+2,2}+u_{i+3,3}+\cdots
+u_{t+i-1,t-1})=$$ $$ (u_{i+1,2}-u_{i+1,1})+
(u_{i+2,3}-u_{i+2,2})+ \cdots +(u_{t+i-1,t}-u_{t+i-1,t-1})\ge 0,$$
and
$$u_{t+c-1,t}\ge u_{t+r,t} ,$$
 we obtain
$$u_{t+c-1,t}\prod_{i=1}^{c-1}m_i\ge (u_{t+r,t}m_1\cdots
m_rm_{r+2}'\cdots m_{c-1}'m_c')
$$
for all $r$, $0\le r \le c-2$, and the lower bound follows.

\vskip 2mm (2) The upper bound is proved similarly. We again
proceed by induction on the codimension $c$ of $I$. If $1\le c \le
2$ then the result works. So, let us assume $c\ge 3$. We will now
induct on $t$. If $t=1$ then $I$ is a complete intersection ideal
and the result is true. Assume $t>1$. Let $\cA''$ (resp. $\cC$) be
the matrix obtained deleting the first column and the first row
(resp. the first column) of the matrix $\cA$ and let $I''$ (resp.
$K$) be the ideal generated by the maximal minors of $\cA''$
(resp. $\cC$). Let $M_i$, $M_i''$ and $M_i(K)$ be the maximal
shifts in the graded minimal free $R$-resolution of $R/I$, $R/I''$
and $R/K$, respectively.

\vskip 2mm By Lemma \ref{key} (3), $$e(R/I)=e(R/I'')+u_{1,1}\cdot
e(R/K),$$ by hypothesis of induction on $c$ and Lemma
\ref{keyinduction} (4), we have
$$e(R/K)\le
\frac{\prod_{i=1}^{c-1}M_i(K)}{(c-1)!}=\frac{\prod_{i=1}^{c-1}M_i}{(c-1)!},$$
and by hypothesis of induction on $t$ we have $$e(R/I'')\le
\frac{\prod_{i=1}^cM_i''}{(c)!}.$$

By Lemma \ref{keyinduction} (2),
$M_i=M_i(I)=M_i(I'')+u_{c-i+1,1}=M_i''+u_{c-i+1,1}$ for all $1\le
i \le c$. Therefore,  we have
$$c!e(R/I)\le \prod_{i=1}^cM_i$$ if and only if
$$\prod_{i=1}^cM_i''+cu_{1,1}\prod_{i=1}^{c-1}M_i\le
\prod_{i=1}^cM_i=$$
$$\prod_{i=1}^c(M_i''+u_{c-i+1,1})=$$

 $$u_{1,1}\prod_{i=1}^{c-1}M_i+
\prod_{i=1}^cM_i''+M_c''\sum_{r=0}^{c-2}(u_{c-r,1}M_1\cdots
M_rM_{r+2}''\cdots M_{c-1}'')=$$
$$u_{1,1}\prod_{i=1}^{c-1}M_i+
\prod_{i=1}^cM_i''+\sum_{r=0}^{c-2}(u_{c-r,1}M_1\cdots
M_rM_{r+2}'\cdots M_{c-1}''M_c'')
 $$ if and only if
$$(c-1)u_{1,1}\prod_{i=1}^{c-1}M_i\le \sum_{r=0}^{c-2}(u_{c-r,1}M_1\cdots
M_rM_{r+2}''\cdots M_{c-1}''M_c'').
$$

Because, for all integer $i$, $ 1\le i \le c-1$, and all integer
$r$, $ 0\le r \le c-2$, we have

$$M_{i}-M_{i+1}''= (u_{c-i+1,1}+u_{c-i+2,2}+\cdots
+u_{t+c-i-1,t-1} +u_{t+c-i,t}+u_{t+c-i+1,t}+\cdots +u_{t+c-1,t})-
$$
$$(u_{c-i+1,2}+u_{c-i+2,3}+\cdots
+u_{t+c-i-1,t} +u_{t+c-i,t}+u_{t+c-i+1,t}+\cdots +u_{t+c-1,t})=$$
$$ (u_{c-i+1,1}-u_{c-i+1,2})+ (u_{c-i+2,2}-u_{c-i+2,3})+
\cdots +(u_{t+c-i-1,t-1}-u_{t+c-i-1,t})\le 0,$$ and
$$u_{1,1}\le u_{c-r,1} ,$$
 we deduce
$$u_{1,1}\prod_{i=1}^{c-1}M_i\le (u_{c-r,1}M_1\cdots
M_rM_{r+2}''\cdots M_{c-1}''M_c'')
$$
for all $r$, $0\le r \le c-2$. This completes the proof of the
upper bound and hence of the Theorem.
\end{proof}

We will end this note with an example which illustrate that the
bounds given in Theorem \ref{maintheorem} are optimal.

\begin{example} Let $I\subset R$ be a codimension $c$
determinantal ideal  generated by the maximal minors of a $t\times
(t+c-1)$ matrix all whose entries are homogeneous polynomials of
fixed degree $1\le d\in \ZZ$. Thus, we have
$$m_i(I)=M_i(I)=td+(i-1)d \quad \text{ for all } i, \quad 1\le i \le c.$$

Therefore, we conclude that
$$
e(R/I)=\frac{\prod_{i=1}^cm_i(I)}{c!}=\frac{\prod_{i=1}^cM_i(I)}{c!}=$$
$$\frac{\prod_{i=1}^c(td+(i-1)d)}{c!}=d^c{t+c-1\choose c}.$$
\end{example}

\end{document}